

\documentstyle{amsppt}
\magnification=1200

\loadbold
\def\Co#1{{\Cal O}_{#1}}
\def\Fi#1{\Phi_{|#1|}}
\def\rw{\rightarrow}
\def\KX{K_X}

\def\Bbbp1{{\Bbb P}^1}
\def\simlin{\sim_{\text{lin}}}
\def\simnum{\sim_{\text{num}}}
\def\Pic{\text{Pic}}

\NoRunningHeads
\TagsOnRight
\pagewidth{14cm}
\pageheight{20cm}
\vcorrection{-0.30in}
\hcorrection{-0.15in}
\topmatter
\title Complex varieties of general type whose canonical systems are
  composed with pencils
\endtitle
\author By Meng Chen$^*$ \endauthor
\thanks $^*$ Supported by NNSFC \endthanks
\address -----------
\newline Department of Applied Mathematics
\newline TongJi University
\newline Shanghai,\ 200092\ \ China \endaddress
\endtopmatter
\document
\baselineskip 12pt
\vskip1cm

Throughout this paper, most our notations and terminologies are standard
within algebraic geometry except for the following which we are in favour
of:

$:=$ --- definition;

$\simlin$ --- linear equivalence;

$\simnum$ --- numerical equivalence.

Let $X$ be a complex nonsingular projective variety of general type with
dimension $d$($d\ge 2$). Suppose  $\text{dim}\Fi{\KX}(X)=1$, we usually say that
the canonical system $|\KX|$ is composed with a pencil. Taking possible 
blow-ups $\pi:X'\rw X$ according to Hironaka such that $g:=\Fi{\KX}\circ\pi $
is a morphism. We have the following commutative diagram:

$$\CD
X'   @>f>>    C\\
@|            @VV{\psi}V\\
X'   @>>g>   W_1\subset{\Bbb P}^{p_g-1}\\
@V{\pi}VV            @.\\
X    @.
\endCD$$
where we set $W_1:=\overline{\Fi{\KX}(X)}$ and let $g:=\psi\circ f$ be the
Stein factorization of $g$. Note that $f$
is a fibration onto a nonsingular curve $C$. 
Let $F$ be a general fiber of $f$, then $F$ is a nonsingular
projective variety of dimension $d-1$. We also say that $f$ is a derived
fibration of the canonical map. Denote $b:=g(C)$, the genus of $C$. 

The aim of this note is to build the following theorems.

\proclaim{\smc Theorem 1} Let $X$ be a complex nonsingular projective
variety of general type with dimension $d$($d\ge 3$). Suppose the canonical
system $|\KX|$ be composed with a pencil, using the above diagram and notations,
and assume that $b\ge 2$, then either
$$p_g(F)=1,\ p_g(X)\ge b-1,$$
or
$$b=p_g(F)=p_g(X)=2.$$
\endproclaim

\newpage
\remark{\smc Remark 1}
In the case of dimension 2, one can refer to \cite{10}.
\endremark

\proclaim{\smc  Theorem 2} 
Under the same assumption as in Theorem 1, assume in addition that 
$\text{dim}X=3$ and $K_X$ is nef and big, set $F_1:={\pi}_*F$,
then

(1) If $b=1$, then $p_g(F)\le 38$; 

(2) If $b=0$, $K_X\cdot F_1^2=0$ and $p_g(X)\ge 20$, then
$$p_g(F)\le 38+\frac{756}{p_g(X)-19};$$

(3) If $b=0$ and $K_X\cdot F_1^2>0$, then
$$p_g(F)\ge q(F)\ge\frac{1}{36}(p_g(X)-1)(p_g(X)-37).$$
\endproclaim
The author would like to thank the referee for many skillful suggestions
which greatly simplify our proofs and improve (1) of Theorem 2 to the present
form. Also, I am indebt to Prof. Zhijie Chen whose opinion improves my calculation
in (2) of Theorem 2.

\head \S1. Proof of Theorem 1\endhead
\demo{Proof of Theorem 1}
Using the first commutative diagram of this paper, we have the fibration
$f:X'\rw C$. Denote by ${\Cal L}$ the saturated subbundle of $f_*\omega_{X'}$
which is generated by $H^0(C,f_*\omega_{X'}),$ ${\Cal L}$ is of rank one under
the assumption of the theorem. Thus we obtain the following exact sequence
$$0\rw {\Cal L}\rw f_*\omega_{X'}\rw {\Cal Q}\rw 0.$$
We know that $R^if_*\omega_{X'/C}$ is a semi-positive vector bundle(Griffith-
Fujita-Kawamata-Ohno semipositivity). ${\Cal Q}\otimes\omega_C^{-1}$ is
automatically semi-positive. Since $H^0(C,{\Cal L})=H^0(C,f_*\omega_{X'}),$
$H^0(C,{\Cal Q})$ injects into $H^1(C,{\Cal L}).$ By Riemann-Roch, 
$h^0(C,{\Cal Q})\ge (b-1)(p_g(F)-1),$ $h^1(C,{\Cal L})\le b-1,$ and hence
we get $p_g(F)\le 2.$ 

If $p_g(F)=1$, then $p_g(X)=h^0(C,{\Cal L})\ge b-1$ by the semipositivity
of ${\Cal L}\otimes\omega_C^{-1}=f_*\omega_{X'/C}.$

If $p_g(F)=2,$ then we must have $h^0({\Cal Q})=h^1({\Cal L})=b-1,$
so that ${\Cal Q}\otimes\omega_C^{-1}$ is of degree $\le 0.$ Since
$f_*\omega_{X'}\otimes\omega_C^{-1}$ is semi-positive, this means that ${\Cal L}$
is of degree $\ge 2b-2,$ with non-trivial $H^1$. Hence ${\Cal L}\cong \omega_C,$
and $1=h^1(C,{\Cal L})\ge h^0({\Cal Q})=b-1.$ This shows that $b=2,$ 
$p_g(X)=h^0(C,{\Cal L})=2,$ completing the proof of Theorem 1.
\enddemo

\head \S2. Proof of Theorem 2\endhead
At first, let us recall Miyaoka's inequality as follows.
\proclaim{\smc Fact } (\cite{7}) Let $X$ be a nonsingular projective
threefold with nef and big canonical divisor $K_X$. Then $3c_2-c_1^2$ is
pseudo-effective, where $c_i$'s are Chern numbers of $X$.
\endproclaim
\proclaim{\smc Proposition 2.1} Let $X$ be a nonsingular projective 3-fold
with nef and big canonical divisor $K_X$. Assume that $|K_X|$ be composed with a 
pencil and that $K_X\cdot F_1^2=0$, then 
$$\Co{F}(\pi^*(K_X)|_F)\cong\Co{F}(\sigma^*(K_{F_0})),$$
where $\sigma:F\rw F_0$ is the contraction onto the minimal model.
\endproclaim
\demo{Proof} This can be obtained by a similar argument as that of Theorem 7 in
\cite{6}.\enddemo

\demo{\smc Proof of Theorem 2} 
We use the same exact sequence as in the proof of Theorem 1. It is easy to
see that $K_{X'}\simnum (\deg{\Cal L})F+Z,$ where $Z$ is the fixed part.

Case (1). $b=1$. In this case, we can suppose $X'=X$. ${\Cal Q}$ is semi-positive.
Note that $\deg {\Cal Q}=0$; otherwise, by Riemann-Roch, we should have
$$h^0(C,f_*\omega_X)\ge \deg {\Cal L}+\deg {\Cal Q}>\deg {\Cal L}=h^0(C,{\Cal L}),$$
a contradiction. $R^1f_*\omega_X$ is semi-positive, while $R^2f_*\omega_X\cong
\Co{C}$ by duality. Thus
$$\align
\chi(X,\omega_X)&=\chi(C,f_*\omega_X)-\chi(C,R^1f_*\omega_X)+\chi(C,R^2f_*\omega_X)\\
&\le\chi(C,f_*\omega_X)=\deg {\Cal L}=p_g(X).
\endalign$$
It follows from Miyaoka's inequality that $K_X^3\le 72\deg {\Cal L}$.
On the other hand,
$$\align
 K_X^3&\ge (\deg {\Cal L}) K_X^2\cdot F\\
&=(\deg {\Cal L}) K_F^2\ge 2(\deg {\Cal L})(p_g(F)-2).
\endalign$$
This implies that $(\deg {\Cal L})(p_g(F)-2)\le 36\deg {\Cal L},$
i.e. $p_g(F)\le 38$.

Case (2). $b=0$ and $K_X\cdot F_1^2=0$. In this case, any vector bundle
on $C$ is a direct sum of  line bundles. $f_*\omega_{X'}$ is a direct sum of ${\Cal L}$
and $p_g(F)-1$ line bundles of degree $-1$ or $-2$, while
$R^1f_*\omega_{X'}$ is a direct sum of $q(F)$ line bundles of degree $\ge -2$.
Hence
$$\align
\chi(X',\omega_{X'})&=\chi(C,f_*\omega_{X'})-\chi(C,R^2f_*\omega_{X'})+
  \chi(C,\omega_C)\\
&\le \deg {\Cal L}+1+q(F)-1.
\endalign$$
Then, by Proposition 2.1, we get
$$(\deg {\Cal L})(p_g(F)-2)\le\frac{1}{2}(\deg {\Cal L})K_{F_0}^2
 \le 36(\deg {\Cal L}+q(F)).$$
Apply the inequality $p_g(F)\ge 2q(F)-4$ (\cite{2}), we get
$$(\deg {\Cal L}-18)p_g(F)\le 38 \deg {\Cal L}-72.$$
We obtained the desired inequality by substituting $\deg {\Cal L}$ by
$p_g(X)-1$. 

Case (3). $b=0$ and $K_X\cdot F_1^2>0$. It is easy to check that
$K_X\cdot F_1^2$ is an even number. We get the following inequality
$$K_X^3\ge (\deg {\Cal L})^2K_X\cdot F_1^2\ge 2(\deg {\Cal L})^2.$$
Therefore we have
$$2(\deg {\Cal L})^2\le 72(\deg {\Cal L}+q(F)),$$
which directly induces what we want. The proof is completed.
\enddemo

\head \S3. Examples\endhead

\subhead 3.1 Examples with $p_g(F)=1$ and $  p_g(X)\ge b$\endsubhead

Let $S$ be a minimal surface with $K_S^2=1$ and $p_g(S)=q(S)=0$.
We  chose $S$ in such a way that
$S$ admits a torsion element $\eta$ of order 2, i.e., $\eta\in
Pic^0 S$ and $2\eta\simlin 0$. For the
existence of this surface, we may refer to \cite{1} or
\cite{9}. 
 Let $C$ be a nonsingular curve of genus $b\ge 1$. Set
$Y:=C\times S$. 
Let $p_1:Y\longrightarrow C$ and $p_2:Y\longrightarrow S$ be the two
projection maps. Let $D$ be an effective divisor on $C$ with
$\deg D=a>0$.  Let $\delta=p_1^*(D)+p_2^*(\eta)$,
$B\simlin
2\delta\simlin p_1^*(2D)$, we can take
$B$ composed of $2a$ distinct points. The pair $(\delta, B)$
determines a smooth double cover $X$ over $Y$. We have the
following commutative diagram, where $\pi$ is the double covering.

$$\CD
X   @>{\pi}>>    {Y=C\times S}@>p_2>> S\\
@|            @VVp_1V    @.\\
X   @>>{\Psi}>   C       @.    @.
\endCD$$
We can see that 
$\Phi_{|K_X|}$ factors through $\Psi$. 
These examples satisfy
$p_g(X)=a+b-1$,  $K_X^3=6(a+2b-2)$,
$h^2({\Cal O}_X)=0$,
$q(X)=b$, $K_F^2=2$, $p_g(F)=1$ and $q(F)=0$.

\subhead 3.2 Examples with $p_g(F)=1$ and $p_g(X)=b-1$\endsubhead

In the construction of 3.1, take $C$ be a hyperelliptic
curve of genus $b\ge 3$. Let $\tau=p-q$ be a divisor on $C$ such
that $2\tau\simlin 0$. Take
$\delta=p_1^*(\tau)+p_2^*(\eta)$, then
$2\delta\simlin 0$. Therefore $\delta$
determines an unramified double cover $\pi:X\longrightarrow Y$.
$X$ is an example with
$p_g(X)=b-1$, $h^2({\Cal O}_X)=0$,
$q(X)=b$ and $K_X^3=6(2b-2)$.  

\subhead 3.3 Examples with $b=0$ and $p_g(F)=2$\endsubhead

In the construction of 3.1,  take $S$ be a minimal surface $S$ with 
$K_S^2=2$ and
$p_g(S)=q(S)=1$ and take $C={\Bbb P}^1$.
Take an effective divisor $D$ on
$C$ with $\deg D=a\ge 3$, $\eta\in Pic^0 S$ with
$2\eta\simlin 0$. Denote by
$\delta:=p_1^*(D)+p_2^*(\eta)$ and $R\simlin
2\delta\simlin p_1^*(2D)$. 
 Thus the pair
$(\delta, R)$ determines a double covering $\pi:X\rightarrow Y$.
We can check that $X$ is an example with
 $p_g(X)=a-1$, $K_X^3=6(a-2)$, $b=0$, $K_F^2=4$, $q(F)=1$ and $p_g(F)=2$.

\subhead 3.3 Examples with $b=1$ and $p_g(F)=2$\endsubhead

Let $S$ be a minimal surface with $p_g(S)=q(S)=1$ and $K_S^2=2$. We know that
the albanese map of $S$ is just a genus two fibration onto an elliptic curve.
It can be constructed from a double cover onto a ruled surface $P$ which
is over the elliptic curve with invariant $e=-1$. Furthermore, all  the
singularities on the branch locus corresponding to this double cover 
are negligible. Let $\pi_0:S\rw \tilde P$ be this double cover with
covering data $(\delta_0,R_0)$. Note that $q(P)=1$.

Take an elliptic curve $E$ and denote $T:=E\times\tilde P$. Let $p_1$ and $p_2$
be two projection maps. Take a 2-torsion element $\eta\in \Pic^0E$.
Since $\tilde P$ is a fibration over an elliptic curve, we can take a 2-torsion 
element $\tau\in \Pic^0\tilde P$ such that $\pi_0^*(\tau)\not\simlin 0$
through the double cover $\pi_0$. 

Let $\delta_1=p_1^*(\eta)+p_2^*(\delta_0)$ and $R_1=2\delta_1$, then the pair
$(\delta_1,R_1)$ determines a double cover $\Pi_1:Y\rw T$. Let 
$\phi=p_1\circ\Pi_1$. Take a divisor $A$ on $E$ with $\deg A=a>0$.
Let $\delta_2=\phi^*(A)+\Pi_1^*p_2^*(\tau)$ and $R_2=2\delta_2$,
then the pair $(\delta_2, R_2)$ determines a smooth double cover 
$\Pi_2:X\rw Y$. We can see that $X$ is a minimal threefold of general type
and the canonical system $|\KX|$ is composed with a pencil. 
$\Fi{\KX}$ factors through $\Pi_2$ and $\phi$. This example satisfies
$b=1$, $p_g(F)=2$, $q(F)=1$, $K_F^2=16$, $p_g(X)=a$, $q(X)=2$
and $K_X^3=12a$.

\head References \endhead
\roster
\item"[1]" Barth, W., Peters, C., Van de Ven, A.: {\it Compact complex
	surfaces,} Berlin Heidelberg New York: Springer 1984.
\item"[2]" Debarre, O.: Addendum to In\'egalit\'es num\'eriques pour
les surfaces de type g\'en\'eral, {\it Bull. Soc. Math. France} {\bf
111}/4, 301-302(1983).
\item"[3]" Hironaka, H.: Resolution of singularities of an algebraic
	variety over a field of characteristic zero, {\it Ann. of
Math.} {\bf 79} (1964), 109--326.
\item"[4]" Kawamata, Y.: Kodaira dimension of algebraic fiber
spaces over curves, {\it Invent. Math.} {\bf 66} (1982), 57--71.
\item"[5]" Koll\'ar, J.: Higher direct images of dualizing
sheaves II, {\it Ann. of Math.} {\bf 124} (1986), 171--202.
\item"[6]" Matsuki, K.: On pluricanonical maps for threefolds of
general type, J. Math. Soc. Japan, Vol. 38, No. 2, 1986, 339-359.
\item"[7]" Miyaoka, Y.: The pseudo-effectivity of $3c_2-c_1^2$
for varieties with numerically effective canonical classes,
In: {\it Algebraic Geometry, Sendai, 1985} (Adv. Stud. in Pure Math.
10, 1987, pp. 449--476).
\item"[8]" Ohno, K.: Some inequalities for minimal fibrations of surfaces
of general type over curves, {\it J. Math. Soc. Japan} {\bf 44} (1992), 643--666.
\item"[9]" Reid, M.: Surfaces with $p_g=0$, $K^2=1$, {\it J. Fac. Sci. Uni.
Tokyo}, Sect. IA, {\bf 25} (1978), 75--92.
\item"[10]" Xiao, G.: L'irr\'{e}gularit\'{e} des surfaces de
type g\'{e}n\'{e}ral dont le syst\`eme canonique est compos\'{e}
d'un pinceau, {\it Compositio Math.} {\bf 56} (1985), 251--257.
\endroster

\enddocument